\documentclass[reqno]{amsart}
\usepackage{amssymb}
\usepackage{amsmath}

\makeatletter
\@addtoreset{equation}{section}
\makeatother

\renewcommand\thefigure{\thesection.\@arabic\c@figure}
\renewcommand\thetable{\thesection.\@arabic\c@table}
 
\newtheorem{theorem}{Theorem}[section]
\newtheorem{lemma}[theorem]{Lemma}
\newtheorem{proposition}[theorem]{Proposition}
\newtheorem{corollary}[theorem]{Corollary}

\newtheorem{remark}[theorem]{Remark}
 
\newcommand{\mc}[1]{{\mathcal #1}}

\newcommand{\bb}[1]{{\mathbb #1}}

\newcommand{\<}{\langle}
\renewcommand{\>}{\rangle}

\def\emptysquare{{\hbox{\vrule height6pt width0.6pt depth0pt%
\vbox{\hrule height0.6pt width4.8pt depth0pt%
\vglue4.8pt%
\hrule height0.6pt width4.8pt depth0pt}%
\vrule height6pt width0.6pt depth0pt}}}
 
\def\qed{\unskip\nobreak
\hfil\penalty50\hskip1.75em\null\nobreak\hfil\emptysquare
{\parfillskip=0pt \finalhyphendemerits=0 \par}\medskip}

\begin{document}
 
\title[Fluctuations in boundary driven exclusion processes]{Stationary
  and nonequilibrium fluctuations in boundary driven exclusion
  processes}

\author{C. Landim, A. Milan\'es, S. Olla}

\address{\noindent IMPA, Estrada Dona Castorina 110, CEP 22460 Rio de
  Janeiro, Brasil and CNRS UMR 6085, Universit\'e de Rouen, UMR 6085,
  Avenue de l'Universit\'e, BP.12, Technop\^ole du Madrillet, F76801
  Saint-\'Etienne-du-Rouvray, France.  
\newline e-mail: \rm
  \texttt{landim@impa.br} } 

\address{Departamento de Estat\'\i stica, ICEx, UFMG,
Campus Pampulha, CEP 31270-901, Belo Horizonte, Brasil
\newline
e-mail: \rm \texttt{aniura@est.ufmg.br}}

\address{Ceremade, UMR CNRS 7534,
  Universit\'e de Paris - Dauphine, Place du Mar\'echal De Lattre
  De Tassigny 75775 Paris Cedex 16 - France.  
\newline e-mail: \rm
  \texttt{olla@ceremade.dauphine.fr}}

\thanks{This research has been partially supported by the agreement
  France-Br\'esil; by the ACI 168 ``Tranport hors \'equilibre'' du
  Minist\`ere de l'\'Education National, France; by FAPERJ Cientistas
  do Nosso Estado; by CNPq Edital Universal and by PRONEX}

\subjclass[2000]{60K35, 82C22}

\keywords{Stationary nonequilibrium states, hydrodynamic fluctuations,
  boundary driven systems, exclusion process}

\begin{abstract} 
  We prove nonequilibrium fluctuations for the boundary driven
  symmetric simple exclusion process. We deduce from this result the
  stationary fluctuations.
\end{abstract}

\date{\today}

\maketitle
 
\section{Introduction}
\label{sec0}     

In the last years there has been considerable progress in
understanding stationary non equilibrium states (SNS): reversible
systems in contact with reservoirs imposing a gradient on the
conserved quantities of the system. In particular, large deviation
properties has been studied for boundary driven one-dimensional
symmetric simple exclusion processes (\cite{bdgjl, delo} and
references therein).

One of the most striking typical property of SNS is the presence of
long-range correlations. For the symmetric simple exclusion this was
already shown by H. Spohn in the pioneering paper \cite{s}. But we
noticed that a mathematical proof of the convergence of the
fluctuation fields to the corresponding Gaussian field was missing
from the literature. The purpose of this paper is to fill this gap.

We consider the symmetric exclusion process in an open lattice of
length N. Particles jumps to nearest neighbors performing simple
symmetric random walks with the exclusion rule: a jump is suppressed
if site is already occupied. At the left boundary particles are
created with rate $\alpha$ and annihilated with rate $1-\alpha$. On
the right boundary this is done with rates $\beta$ and $1-\beta$. To
keep notation simple we restrict to the one dimensional nearest
neighbor case. Extensions to more dimension and more general jumps
rates is straightforward (see remarks \ref{so1} and \ref{s02}).

If $\alpha = \beta = \rho$, the Bernoulli product measure with
probability $\rho$ is stationary and reversible for the dynamics. But
when $\alpha \neq \beta$ the stationary measure has correlations
and is not explicitly computable. We denote by $\eta(x) = 0$ or $1$
the occupation variable of site $x$. It is easy to prove that
$<\eta([Nu])>_{ss} \to \bar\rho(u) = (\beta - \alpha) u + \alpha$.

The fluctuation field is formally defined as the random distribution
on $[0,1]$
\begin{equation}
  \label{eq:1}
  Y^N(u) = \frac 1{\sqrt N} \sum_{x=1}^{N-1} \delta (u - x/N) (\eta(x) -
  \bar\rho(u))\; . 
\end{equation}
We prove in this paper that, under the stationary measure for the
process, $Y^N$ converges in law to the centered Gaussian field $Y$ on
$[0,1]$ with covariance
\begin{equation}
  \label{eq:2}
  <Y(u) Y(v)> \;=\; \chi(\bar\rho(u)) \delta(u-v) \;-\; 
(\beta - \alpha)^2 (-\Delta)^{-1}(u,v) \;,
\end{equation}
where $\Delta$ is the Laplace operator with Dirichlet boundary
conditions, and $\chi(\rho) = \rho(1-\rho)$.  In the 1-dimensional
case we have more explicitly $(-\Delta)^{-1}(u,v) = u(1-v)$.

The strategy we use to prove this result is to study first the
convergence of the nonstationary fluctuations. If we start with some
non-equilibrium density profile $<\eta_0([N u])> = \rho(0, u)$, then
at the diffusive time scale we have $<\eta_{N^2 t}([N u]) > \to
\rho(t,u)$, where $\rho(t,u)$ is solution of the heat equation with
initial condition $\rho(0,u)$.

We then consider the time-dependent fluctuation field
\begin{equation}
  \label{eq:3}
    Y^N(t, u) = \frac 1{\sqrt N} \sum_{x=1}^{N-1} 
\delta (u - x/N) (\eta_{N^2 t}(x) - \rho(t,u )). 
\end{equation}
The main point of the proof is to show the convergence of $Y^N(u,t)$ to
the solution of the stochastic linear partial differential equation
\begin{equation}
  \label{eq:4}
  \partial_t Y(t,u) = \Delta Y(t,u) - \nabla
  \left(\sqrt{2\chi(\rho(t, u))}\; W(t,u) \right) \;,
\end{equation}
where $W(t,u)$ is the standard space-time white noise.  If we start in
the stationary state, $\rho(t,u) =\bar \rho(u)$ for all $t$.  In this case
the distribution valued process $Y(t,u)$ is a stationary Gaussian
process and its invariant distribution is given by the Gaussian field
$Y$ with covariance given by (\ref{eq:2}).

This article presents a rigorous proof of the results described above
and presented in \cite{delo}. Article \cite{delo} also contains 
the connection between the large deviations and the small fluctuations
proved here, showing
that the inverse of the covariance \eqref{eq:2} is given by the second
functional derivative of the large deviations rate function. 

\section{Notation and results}
\label{sec1}

For $N\ge 1$, let $\Lambda_N =\{1, \dots, N-1\}$.  Fix $0\le \alpha\le
\beta\le 1$ and consider the boundary driven symmetric simple
exclusion process associated to $\alpha$, $\beta$. This is the Markov
process on $\{0,1\}^{\Lambda_N}$ whose generator $L_N$ is given by
\begin{eqnarray*}
(L_N f)(\eta) &=& \sum_{x=1}^{N-2} \{ f(\sigma^{x,x+1}
\eta)-f(\eta)\} \\
&+& \Big\{ \alpha [1-\eta(1)] + (1-\alpha) \eta(1) \Big\}
\{ f(\sigma^{1} \eta)-f(\eta)\} \\
&+& \Big\{ \beta [1-\eta(N-1)] + (1-\beta) \eta(N-1) \Big\}
\{ f(\sigma^{N-1} \eta)-f(\eta)\}\;.
\end{eqnarray*}
In this formula, $\eta = \{\eta(x),\, x\in\Lambda_N\}$ is a
configuration of the state space $\{0,1\}^{\Lambda_N}$ so that $\eta
(x)=0$ if and only if site $x$ is vacant for $\eta$;
$\sigma^{x,y}\eta$ is the configuration obtained from $\eta$ by
interchanging the occupation variables $\eta(x)$, $\eta(y)$:
$$
(\sigma^{x,y} \eta) (z)\; =\;
\left\{
\begin{array}{ll}
\eta (z)  & \hbox{if $z\neq x,y$}\; , \\
\eta (y)  & \hbox{if $ z=x$}\; , \\
\eta (x)  & \hbox{if $ z=y$}\; ;
\end{array}
\right.
$$
and $\sigma^{x}\eta$ is the configuration obtained from $\eta$ by
flipping the variable $\eta(x)$:
$$
(\sigma^{x} \eta) (z)\; =\;
\left\{
\begin{array}{ll}
\eta (z)  & \hbox{if $z\neq x$}\; , \\
1-\eta (z)  & \hbox{if $ z=x$}\; .
\end{array}
\right.
$$
Hence, at rate $\alpha$ (resp. $1-\alpha$) a particle is created
(resp. removed) at the boundary site $1$ if this site is vacant (resp.
occupied). The same phenomenon occurs at the boundary $x=N-1$ with
$\beta$ in place of $\alpha$.

This finite state Markov process is irreducible and has therefore a
unique stationary measure, denoted by $\nu^N_{\alpha, \beta}$. For $0\le
\gamma\le 1$, denote by $\nu^N_{\gamma}$ the Bernoulli product measure
on $\{0,1\}^{\Lambda_N}$ with density $\gamma$.  If $\alpha=\beta$,
an elementary computation shows that $\nu^N_{\alpha}$ is the invariant
measure and that the process is reversible with respect to this
stationary state.  On the other hand, if $\alpha<\beta$, it is known
since \cite{s} that the invariant state has long range correlations.

\medskip
\noindent{\bf Static picture.} For $N\ge 1$, denote by $\pi^N$ the
measure on $[0,1]$ obtained by assigning mass $N^{-1}$ to each
particle:
\begin{equation}
\label{f23}
\pi^N(\eta) \;=\; N^{-1} \sum_{x\in\Lambda_N} \eta(x) \delta_{x/N}\;,
\end{equation}
where $\delta_u$ is the Dirac measure on $u$. It has been proved in 
\cite{els} that under the stationary state $\nu^N_{\alpha, \beta}$ 
the empirical measure $\pi^N$ converges to the unique solution of the
elliptic equation
$$
\left\{
\begin{array}{l}
\Delta \rho = 0\;, \\
\rho(0) = \alpha\;, \quad \rho(1) = \beta \;.
\end{array}
\right.
$$
We denote the solution of this equation by $\bar\rho =
\bar\rho_{\alpha, \beta}$.
\medskip

Once a law of large number has been proved for the empirical measure
under the stationary state, it is natural to consider the fluctuations
around the limit. Let 
$$
\rho^N(x) = E_{\nu^N_{\alpha, \beta}}[\eta(x)]\;.
$$
Since $E_{\nu^N_{\alpha, \beta}}[L_N \eta(x)] =0$ for all $1\le
x\le N-1$, an elementary computation shows that $\rho^N$ is the
solution of
\begin{equation}
\label{eq:04}
\left\{
\begin{array}{l}
(\Delta_N \rho^N)(x) = 0 \quad \text{ for $1\le x\le N-1$}\;, \\
\rho^N(0) = \alpha \; , \quad \rho^N( N ) = \beta \;,
\end{array}
\right.
\end{equation}
where $\Delta_N$ is the discrete Laplacian:
$(\Delta_N H)(x) = N^2\{ H(x+1) + H(x-1) - 2 H(x)\}$. 

In the case of the symmetric simple exclusion process,
$\rho^N(\cdot)$ is just the linear interpolation between
$\rho^N(0)=\alpha$, $\rho^N(N)=\beta$. 

To define the space in which the fluctuations take place, denote by
$C^2_0([0,1])$ the space of twice continuously differentiable
functions on $(0,1)$ which are continuous on $[0,1]$ and which vanish
at the boundary.  Let $- \Delta$ be the positive operator, essentially
self-adjoint on $L^2[0,1]$, defined by
\begin{eqnarray*}
- \Delta&=&-\frac{d^2}{dx^2}\;, \\ 
\mathcal{D}(-\Delta)&=& C^2_0([0,1]) \;. 
\end{eqnarray*}
Its eigenvalues and corresponding (normalized) eigenfunctions have the
form $\lambda_n=(n\pi)^2$ and $e_n(u)=\sqrt{2}\sin(n\pi u)$
respectively, for any $n\in\mathbb{N}$. By the Sturm-Liouville theory,
$\{e_n,\;n\in \mathbb{N}\}$ forms an orthonormal basis of $L^2[0,1]$.

We denote with the same symbol the closure of $-\Delta$ in $L^2[0,1]$.
For any nonnegative integer $k$, we define the Hilbert spaces
$\mathcal{H}_k=\mathcal{D}(\{-\Delta\}^{k/2})$, with inner product
$(f,g)_k=(\{-\Delta\}^{k/2}f$, $\{-\Delta\}^{k/2}g)$, where $(\cdot
,\cdot )$ is the inner product in $L^2[0,1]$. By the spectral
theorem for self-adjoint operators,
\begin{equation*}
\mathcal{H}_k \;=\; \{f\in L^2[0,1]:\;\sum_{n=1}^{+\infty}n^{2k}
(f,e_n)^2<\infty\}\; ,
\end{equation*}
\begin{equation*}
(f,g)_k \;=\; \sum_{n=1}^{+\infty}(n\pi)^{2k}(f,e_n)(g,e_n)\;.
\end{equation*}

Moreover, if $\mathcal{H}_{-k}$ denotes the topological dual space of
$\mathcal{H}_k$,
\begin{eqnarray*}
\mathcal{H}_{-k} \;=\; \{f\in\mathcal{D}'(0,1):\;\sum_{n=1}^{+\infty}n^{-2k}
\<f,e_n\>^2<\infty\},\\ 
(f,g)_{-k} \;=\; \sum_{n=1}^{+\infty}(n\pi)^{-2k} \<f,e_n\> \<g,e_n\>,
\end{eqnarray*}
where $\<f,\cdot\>$ represents the action of the distribution $f$ over
$[0,1]$ on test functions.  



Fix $k>5/2$ and define the density field $Y^N$ on $\mc H_{-k}$ by
\begin{equation}
\label{eq:08}
Y^N (H) \;=\; N^{-1/2} \sum_{x\in\Lambda_N} H(x/N)
\{ \eta(x) - \rho^N(x)\}\;.
\end{equation}

For $k\ge 1$, denote by $q$ the Gaussian probability measure on
$\mc H_{-k}$ with zero mean and covariance given by
\begin{eqnarray}
\label{eq:07}
\!\!\!\!\!\!\!\!\!\!\!\!\!\! &&
E_{q} [Y(H) Y(G)] \;=\; \\
\!\!\!\!\!\!\!\!\!\!\!\!\!\! && \quad
\int_0^1 du\, \chi(\bar \rho(u)) 
\, H(u) \, G(u) \;-\; (\beta - \alpha)^2  \int_0^1 du\,
[(-\Delta)^{-1} H (u)]  \, G(u) \;.
\nonumber
\end{eqnarray}

\begin{theorem}
\label{mt1}
Fix $k>5/2$ and denote by $q_N$ the probability measure on $\mc
H_{-k}$ induced by the density field $Y^N$ defined in \eqref{eq:08}
and the stationary measure $\nu^N_{\alpha, \beta}$. As $N\uparrow
\infty$, $q_N$ converges to $q$.
\end{theorem}

This result follows from Proposition \ref{s5} which is proved in
Section \ref{sec2}.

\begin{remark}
\label{so1}
The same statement holds in higher dimensions for a symmetric
exclusion on the set $\Lambda_N \times \bb T_N^{d-1}$, where $\bb T_N$
is the discrete torus of length $N$. The dynamics is periodic in the
$d-1$ directions orthogonal to the gradient of the density. In this
case the space correlations are given by \eqref{eq:07}, where $\Delta$
is the Laplacian with periodic boundary conditions in the (d-1)
dimensions and Dirichlet boundary conditions in the first coordinate.
The proof is an elementary extension of the one-dimensional case.
\end{remark}

\begin{remark}
\label{s02}
We may also consider a boundary driven symmetric simple exclusion
process in which the occupation variables $\eta(x)$, $\eta(x+y)$ are
exchange at rate $p(y)$ for a finite range irreducible probability
$p(\cdot)$. In this case, the Laplacian is replaced by the operator
$\sum_{i,j=1}^d \sigma_{i,j} \partial_{u_i} \partial_{u_j}$, where
$\sigma_{i,j} = \sum_y y_i y_j p(y)$.
\end{remark}

\section{Nonequilibrium fluctuations}
\label{sec2}   

We prove in this section the dynamical nonequilibrium fluctuations of
the boundary driven exclusion process. We start with the law of large
numbers. 

Fix a density profile $\rho_0 : [0,1] \to [0,1]$. Consider a sequence
$\{\mu^N,\, N\ge 1\}$ of probability measures on $\{0,1\}^{\Lambda_N}$
such that for every continuous test function $H :[0,1]\to \bb R$ and
every $\varepsilon>0$,
$$
\lim_{N\to\infty}
\mu^N \Big\{ \pi^N( H) - \int H(u) \rho_0(u) du \Big|
> \varepsilon \Big\} \;=\; 0\;.
$$

Denote by $\bb P_{\mu^N}$ the probability on the path space
$D(\bb R_+, \{0,1\}^{\Lambda_N})$ induced by the Markov process with
generator $L_N$ and the initial measure
$\mu^N$. Denote by $\pi^N_t$ the empirical measure associated to the
state of the process at time $t$: $\pi^N_t = N^{-1} \sum_x\eta_{N^2t}(x)
\delta_{x/N}$. 
It follows from the usual hydrodynamic limits techniques, 
adapted to the boundary driven context (cf. sections 4 and 5 in
\cite{kl}, and \cite{lms}) that for every $t \ge 0$, every
continuous test function $H :[0,1]\to \bb R$ and every
$\varepsilon>0$,
$$
\lim_{N\to\infty}
\bb P_{\mu^N} \Big\{\Big| \pi^N_t( H) - \int H(u) \rho(t,u) du \Big|
> \varepsilon \Big\} \;=\; 0\;,
$$
where $\rho(t,u)$ is the unique solution of the heat equation
\begin{equation}
\label{f24}
\left\{
\begin{array}{l}
\partial_t \rho = \Delta \rho\;, \\
\rho(0,\cdot) = \rho_0(\cdot) \;, \\
\rho(\cdot , 0) = \alpha\;, \quad  \rho(\cdot , 1) = \beta\;.
\end{array}
\right.
\end{equation}
Furthermore is valid the following \emph{replacement lemma}:
\begin{lemma}\label{rl}
  Let $\Psi(\eta)$ a local function, and $\tilde\Psi(\rho) =
  E_{\nu_\rho}(\Psi)$, where $\nu_\rho$ is the Bernoulli measure with
  probability $\rho$. Let $G(s,u)$ a continuous function on
  $\bb R_+\times [0,1]$. Then 
  \begin{equation*}
     \limsup_{N\to \infty} E\left(\int_0^T ds \left| \frac 1N
        \sum_{x=1}^{N-1} G(s, \frac xN) \Psi(\tau_x \eta_{N^2s}) - \int_0^1
        G(s,u) \tilde\Psi(\rho(s,u)) \; du\right| \right) = 0
  \end{equation*}
\end{lemma}
The proof is given in chapter 5 of \cite{kl}, adapted to the open
boundary situation.
\medskip

We now turn to the fluctuations. Consider a sequence $\{\mu^N : \,
N\ge 1\}$ of probability measures on $\{0,1\}^{\Lambda_N}$. Let
\begin{equation*}
\rho^N(x) \;=\; E_{\mu^N}[\eta(x)]\;,\quad
\varphi^N(x,y) \;=\; E_{\mu^N}[\eta(x); \eta(y)]\; ,
\end{equation*}
for $x$, $y$ in $\Lambda_N$, $x<y$. In this formula, $E_{\mu}[f ; g]$
stands for the covariance of $f$ and $g$: $E_{\mu}[f ; g] = E_{\mu}[f
g] - E_{\mu}[f]E_{\mu}[g]$. We extend the definition of $\rho^N$ and
$\varphi^N$ to the boundary of $\Lambda_N$ by setting 
\begin{equation*}
\rho^N (0) \;=\; \alpha\;, \quad \rho^N (N) \;=\; \beta\;, \quad 
\varphi^N(x,y) \;=\; 0 
\end{equation*}
if $x$ or $y$ does not belong to $\Lambda_N$.  Assume that there
exists a finite constant $C_0$ such that
\begin{equation}
\label{f21}
\sup_{0\le x\le N-1} N |\rho^N(x+1) - \rho^N(x) | \;\le\; C_0\;,\quad
N \max_{\substack{x,y\in\Lambda_N \\ x<y}} |\varphi^N(x,y)| \;\le\; C_0\;.
\end{equation}
Assume furthermore that $\rho^N$ converges weakly to a profile
$\rho_0$ in the sense that for every continuous function
$H:[0,1]\to\bb R$,
\begin{equation}
\label{f22}
\lim_{N\to\infty} \frac 1N \sum_{x\in\Lambda_N} H(x/N) \rho^N(x) \;=\;
\int_0^1 du\, H(u) \rho_0(u)\;.
\end{equation}

It follows from assumptions \eqref{f21}, \eqref{f22} and from
Chebyshev inequality that under $\mu^N$ the empirical measure $\pi^N$,
defined in \eqref{f23}, converges to $\rho_0(u) du$: For every
$\delta>0$ and every continuous function $H:[0,1]\to\bb R$,
\begin{equation*}
\lim_{N\to\infty} \mu^N \Big\{ \, \Big|  \pi^N ( H) - \int_0^1 du\, 
H(u) \rho_0(u) \Big| > \delta \Big\} \;=\; 0\;.
\end{equation*}
In particular, by the law of large numbers stated in the beginning of
this section, for every $t>0$, the empirical measure $\pi^N_t$
converges to the absolutely continuous measure whose density is the
solution of the heat equation \eqref{f24}.

We prove in \eqref{f25} that the stationary state $\nu^N_{\alpha,
  \beta}$ satisfies the assumptions \eqref{f21}, \eqref{f22}. It also
easy to verify that this property is shared by product measures
associated to Lipschitz profiles.

Let $\rho^N_t$ be the solution of the semidiscrete heat equation
\begin{equation}
\label{f26}
\left\{
\begin{array}{l}
\vphantom{\Big\{}
\partial_s \rho^N_s (x) = \Delta_N \rho_s(x) \;, \;\; x\in
\Lambda_N \;; \\
\vphantom{\Big\{}
\rho^N_0(x) = \rho^N (x) \;, \;\; x\in \Lambda_N \;; \\
\vphantom{\Big\{}
\rho^N_s (0) = \alpha \;, \;\; \rho^N_s(N) = \beta\;, 
\;\;\; s\ge 0\;.
\end{array}
\right.
\end{equation}
Fix $k>5/2$ and denote by $Y^N_t$ the density fluctuation field which
acts on smooth functions $H$ in $\mc H_k$ as
$$
Y^N_t (H) \;=\; N^{-1/2} \sum_{x\in\Lambda_N} H(x/N)
\{ \eta_{tN^2}(x) - \rho_t^N(x)\}\;.
$$
Notice that time has been speeded up by $N^2$.  Denote by $Q_N$ the
probability measure on $D([0,T], \mc H_{-k})$ induced by the density
fluctuation field $Y^N$ introduced above and the probability measure
$\mu^N$.

Assumptions \eqref{f21}, \eqref{f22} ensure tightness of the sequence
$Q_N$ and permit to describe the asymptotic evolution of the field $Y$
as the sum of two uncorrelated pieces: a deterministic part
characterized by the heat kernel and a martingale. This is the content
of the first result. Denote by $\{T_s: s\ge 0\}$ the semigroup
associated to the operator $\Delta$.

\begin{proposition}
\label{s3}
Fix $T>0$ and a positive integer $k> 5/2$. The sequence $Q_N$ is tight
on $D([0,T], \mc H_{-k})$ with respect to to the uniform topology.  All
limit points $Q^*$ are concentrated on paths $Y_t$ such that
\begin{equation}
\label{f41}
Y_t (H) \;=\; Y_0 (T_t H) \;+\; W_t(H)\;,
\end{equation}
where $W_t(H)$ is a zero-mean Gaussian variable with variance given by
\begin{equation*}
2 \int_0^t ds\, \int_0^1 du\;  \chi(\rho(s,u)) (T_{t-s} \nabla H)^2(u)\;,
\end{equation*}
and $\rho_s$ is the solution of the heat equation \eqref{f24}.
Moreover, $Y_0$ and $W_t$ are uncorrelated in the sense that
$E_{Q^*}[Y_0 (H) W_t(G)]=0$ for all functions $H$, $G$ in
$C^2_0([0,1])$ and all $0\le t\le T$.
\end{proposition}

\begin{proof}
The proof of tightness of the sequence $Q_N$ is left to the end of
this section. To check the properties of the limit points, fix a
smooth function $H$ in $C^2_0([0,1])$.  An elementary computation
shows that $$
(\partial_t + N^2 L_N) Y^N_t (H) = Y^N_t (\Delta_N H).
$$
Observe that no boundary term appears in the right hand side of the
above equation. In particular, defining
\begin{equation}
\label{eq:09}
\Gamma^{N}_s (H) \;=\; N^2 \big\{ L_N Y^N_s (H)^2 - 2 
Y^N_s (H) L_N Y^N_s (H) \big\}\;,
\end{equation}
it follows that
\begin{eqnarray}
\label{eq:01}
\!\!\!\!\!\!\!\!\!\!\!\!\! &&
M^{1,N}_t (H) \;=\; Y^N_t (H) \;-\; Y^N_0 (H) \;-\; \int_0^t ds
\, Y^N_s (\Delta_N H)  \;, \\
\!\!\!\!\!\!\!\!\!\!\!\!\! && \quad
M^{2,N}_t (H) \;=\; \big\{M^{1,N}_t (H)\big\}^2
\;-\;  \int_0^t ds \, \Gamma^{N}_s (H)\;, \nonumber
\end{eqnarray}
are martingales.  A simple computation shows that
\begin{eqnarray}
\label{eq:02}
\Gamma^{N}_s (H) &=&  
\frac 1N \sum_{x=1}^{N-2} [\eta_{N^2 s}(x+1) - \eta_{N^2 s}(x)]^2
\big\{ (\nabla_N H)(x/N) \big\}^2 \\
\!\!\!\!\!\!\!\!\!\!\!\!\! 
&+& N^{-1} (\nabla_N H) (0)^2 \{\eta_{N^2 s}(1) - \alpha\}^2 \nonumber \\
\!\!\!\!\!\!\!\!\!\!\!\!\! 
&+& N^{-1} (\nabla_N H) ((N-1)/N)^2 \{\eta_{N^2 s}(N-1) - \beta\}^2  \;,
\nonumber
\end{eqnarray}
where $\nabla_N H (x/N)$ stands for the discrete derivative:
$(\nabla_N H)(u) = N\{ H(u+N^{-1}) - H(u)\}$. Observe that the last
two boundary term on the above equations are of order $N^{-1}$.  By
Lemma \ref{rl}, as $N\uparrow\infty$, $\Gamma^{N}_s (H)$ converges to
$ 2 \int_0^1 \chi(\rho(s,u)) (\nabla H(u))^2\; du$, where $\rho$ is
the solution of the heat equation \eqref{f24}.
  
Fix a limit point $Q^*$ of the sequence $Q_N$. It follows from
\eqref{eq:01} that under $Q^*$ for each $H$ in $C^2_0([0,1])$
\begin{equation}
\label{eq:06}
M_t (H) \;=\; Y_t (H)\;-\; Y_0 (H) \;-\; \int_0^t ds
\, Y_s (\Delta  H)  
\end{equation}
is a martingale with deterministic quadratic variation given by
\begin{equation*}
\< M (H)\>_t \;=\; 2 \int_0^t ds  \int_0^1 du\, \chi(\rho(s,u)) 
(\nabla H)^2\; .
\end{equation*}
In particular, for each $H$, $M_t (H)$ is Brownian motion changed in
time.

Consider the semi-martingale $Y_s(T_{t-s} H)$ for $0\le s\le t$. Apply
Ito's formula to derive equation \eqref{f41}, with
\begin{equation*}
W_t(H) \;=\; \int_0^t dM_s(T_{t-s}H)\;.
\end{equation*}
$W_t(H)$ has a Gaussian distribution because the martingales $M_t(H)$
are Gaussian, being a deterministic time-change of a Brownian motion.
The expression for the variance of $W_t(H)$ follows from an elementary
computation, as well as the fact that $W_t(H)$ and $Y_0(G)$ are
uncorrelated. This concludes the proof of the proposition.
\end{proof}

In view of \eqref{f41}, to prove that $Y^N$ converges it remains to
guarantee the convergence at the initial time:

\begin{proposition}
\label{s1} 
Assume that $Y^N_0$ converges to a zero-mean Gaussian field $Y$ with
covariance denoted by $\ll \cdot, \cdot \gg$:
\begin{equation*}
\lim_{N\to\infty} E_{\mu^N} [Y(H) Y(G)] \;=\;
E[Y(H) Y(G)] \;=:\; \ll H, G \gg\; .
\end{equation*}
Then, $Q^N$ converges to a generalized Ornstein-Uhlenbeck process with
covariances given by
\begin{equation}
\label{f11}
\begin{split}
  E[Y_t(H) Y_s(G)] \;=\;& \ll T_t H , T_s G\gg \\
  &\;+\; 2\int_0^s dr\,
  \int_0^1 \chi(\rho(r,u)) \, (\nabla T_{t-r} H)(u) \, (\nabla T_{s-r}
  G)(u)
\end{split}
\end{equation}
for all $0\le s\le t\le T$, $H$, $G$ in $C^2_0([0,1])$.
\end{proposition}

\begin{proof}
  By Proposition \ref{s3}, the sequence $Q^N$ is tight and all limit
  points satisfy \eqref{f41}. Since $W_t$ and $Y_0$ are zero-mean
  Gaussian random variables, so is $Y_t$. To compute the covariance,
  it is enough to remind that the variables are uncorrelated. The first
  piece in formula \eqref{f11} accounts for the covariance between
  $Y_0(T_t H)$ and $Y_0(T_sG)$, while the last one for the covariance
  between $W_t(H)$ and $W_s(G)$.
\end{proof}

A nonequilibrium central limit theorem for the density field follows
from the previous two results for processes starting from local Gibbs
states. Indeed, fix a Lipschitz profile $\gamma:[0,1] \to [0,1]$ such
that $\gamma (0) =\alpha$, $\gamma (1) = \beta$ and denote by
$\nu^N_{\gamma(\cdot)}$ the product measure on $\{0,1\}^{\Lambda_N}$
associated to $\gamma$ so that
\begin{equation*}
\nu^N_{\gamma(\cdot)} \{ \eta (x) = 1\} \;=\; \gamma(x/N)
\end{equation*}
for $x$ in $\Lambda_N$. In this case $\rho^N(x) = \gamma(x/N)$,
$\varphi^N(x,y) =0$ and $\rho_0=\gamma$. The first hypothesis in
\eqref{f21} is satisfied because we assumed $\gamma$ to be Lipschitz.
On the other hand, computing the characteristic functions of $Y^N_0$,
it is easy to show (cf. \cite{kl}) that $Y^N_0$ converges to a
zero-mean Gaussian field with covariance given by
\begin{equation*}
E[Y(H) Y(G)] \;=\; \int_0^1 \chi(\gamma(u)) \, H(u)\, G(u)\; du.
\end{equation*}
Therefore, by Propositions \ref{s3} and \ref{s1} the density field
converges to a generalized Ornstein-Uhlenbeck process: 

\begin{corollary}
\label{s4}
Fix $T>0$ and a positive integer $k> 5/2$. Denote by $Q_N$ the
probability measure on $D([0,T], \mc H_{-k})$ induced by the density
fluctuation field $Y^N$ and the probability measure
$\nu^N_{\gamma(\cdot)}$. Then, $Q^N$ converges to the centered
Gaussian probability measure $Q$ with covariances given by
\begin{equation*}
  \begin{split}
    E[Y_t(H) Y_s(G)] \;=\;& \int_0^1\; du\; \chi(\gamma(u)) \, T_t H(u)\, T_s
    G(u) \\
    & \;+\; 
    2\int_0^s dr\, \int_0^1 \; du\; \chi(\rho(r,u)) \, (\nabla T_{t-r}
    H)(u) \, (\nabla T_{s-r} G)(u) 
  \end{split}
\end{equation*}
for all $0\le s\le t\le T$, $H$, $G$ in $C^2_0([0,1])$, where $\rho_t$
is the solution of the heat equation with initial condition $\gamma$.
\end{corollary}

A similar result was obtained by De Masi et al. \cite {de}, for the
one-dimensional symmetric exclusion process in infinite volume. This
result was extended to higher dimensions by Ravishankar \cite{r}. Chang
and Yau \cite{cy} introduced a general method to prove non-stationary
fluctuations of one-dimensional interacting particle systems.

In Proposition \ref{s1}, the asymptotic behavior of the covariance
\eqref{f11} as $t\uparrow \infty$ can be computed. Indeed, fix a
function $H$ in $C^2_0([0,1])$ and set $G=H$, $s=t$. Since $T_t H$
vanishes as $t\uparrow \infty$, the first part of the covariance
converges to $0$. Since $H$ vanishes at the boundary, since $T_t$ is
the semigroup associated to the Laplacian and since $2 H \nabla H
= \nabla H^2$, an integration by parts shows that the second part
of the covariance \eqref{f11} is equal to
\begin{equation*}
  \begin{split}
    2 \int_0^t ds \, \, & \int_0^1 \; du\; \chi(\rho(s,u)) \, (T_{t-s} H)(u)
    \, (\partial_s T_{t-s} H)(u) \\
    & \qquad \qquad \;-\; \int_0^t ds \, \, \int_0^1 \; du\;
    [\nabla \chi(\rho(s,u))] \, [\nabla (T_{t-s} H(u))^2 ] \;.
  \end{split}
\end{equation*}
Since $2 G \partial_s G = \partial_s G^2$, integrating by parts in
time, since $T_t H$ vanishes in the limit $t\uparrow\infty$ and since
the solution of the heat equation converges to the stationary profile
$\bar\rho$, the previous expression is equal to
\begin{equation*}
  \begin{split}
    \int_0^1 \; du\; \chi(\bar \rho(u)) H(u)^2 \;-\; \int_0^t ds \,
    \int_0^1 \; du\; [\partial_s \chi(\rho(s,u))] \, (T_{t-s} H(u))^2
    \\
    \; +\; \int_0^t ds \, \int_0^1 \; du\; [\Delta \chi(\rho(s,u))] \,
    (T_{t-s} H(u))^2
  \end{split}
\end{equation*}
plus a term which vanishes in the limit. This sum is equal to
\begin{equation*}
\int_0^1 \; du\;  \chi(\bar \rho(u)) H(u)^2   \;- 2 \; \int_0^t ds \, 
\int_0^1 \; du\; [\nabla \rho(s,u)]^2 \, (T_{t-s} H(u))^2  
\end{equation*}
because $\rho$ is the solution of the heat equation. As
$t\uparrow\infty$, this expression converges to
\begin{equation*}
\int_0^1 \; du\; \chi(\bar \rho(u)) H^2(u)  \;- (\beta-\alpha)^2 \,
\int_0^1 \; du\; H(u)  \, (-\Delta)^{-1} H (u)  
\end{equation*}
because $\nabla \bar\rho=\beta-\alpha$. We just recovered the
covariance \eqref{eq:07} of the density field under the stationary
state. 

We turn now to the proof of Theorem \ref{mt1}.  Assume that the
initial state $\mu^N$ is the stationary state $\nu^N_{\alpha, \beta}$.
We prove in Section \ref{sec3} that the second condition in
\eqref{f21} is fulfilled. 

Fix $k$ in $\bb R$ and recall the definition of the probability
measure $q$ introduced just before \eqref{eq:07}. Let $Q$ be the
probability measure on $C([0,T], \mc H_{-{k}})$ corresponding to the
stationary generalized Ornstein--Uhlenbeck process with mean $0$ and
covariance given by
\begin{equation}
\label{f10}
E_{Q}\Big[ Y_t(H) Y_s(G)\Big]\; = \; E_{q}\Big[ Y(T_{t-s} H) Y(G)\Big]
\end{equation}
for every $0\le s\le t$ and $H$, $G$ in $\mc H_{k}$.

\begin{proposition}
\label{s5}
Fix $T>0$ and a positive integer $k> 5/2$. Denote by $Q_N$ the
probability measure on $D([0,T], \mc H_{-k})$ induced by the density
fluctuation field $Y^N$ and the probability measure $\nu^N_{\alpha,
  \beta}$. The sequence $Q_N$ converges weakly to the probability
measure $Q$.
\end{proposition}

\begin{proof}
  Since $\nu^N_{\alpha, \beta}$ satisfies assumptions \eqref{f21},
  \eqref{f22}, by Proposition \ref{s3}, $Q^N$ is tight and all limit
  points satisfy \eqref{f41}, where $W_t(H)$ is in this stationary
  context a zero-mean Gaussian variable with variance given by
  \begin{equation*}
    2 \int_0^t ds\, \int_0^1 \; du\; \chi(\bar \rho(u)) (T_s \nabla
    H(u))^2\; .
  \end{equation*}
  As $t\uparrow\infty$, $T_t H$ vanishes in $L^2([0,1])$. On the other
  hand, the computations performed just before the statement of this
  lemma show that the variance of $W_t(H)$ converges to \eqref{eq:07}.
  Therefore, $W_t(H)$ converges in distribution to a zero-mean
  Gaussian variable with variance given by \eqref{eq:07}.  Since the
  process is stationary, we just proved that the variables $\{ Y_t(H):
  t\ge 0, H\in\mc H_k\}$ have a zero mean Gaussian distribution with
  covariance given by \eqref{eq:07}.

  To compute the covariances $E_{Q^*}[Y_t(H) Y_s(G)] =
  E_{Q^*}[Y_{t-s}(H) Y_0(G)]$, $0\le s\le t$ it is enough to iterate
  relation \eqref{eq:06} to recover formula \eqref{f10}. This concludes
  the proof of the lemma.
\end{proof}

\medskip
We conclude this section proving that the sequence of probability
measures $Q_N$ is tight and that all limit points are concentrated on
continuous paths.

To prove that the sequence $Q_N$ is tight we need to show that
for every $0\le t\le T$,
$$
\lim_{A\to\infty} \limsup_{N\to\infty} 
\bb P_{\mu^N} \Big[ \sup_{0\le t\le T} \Vert
Y_t\Vert_{-k} >A\Big]\; =\; 0
$$
and that
$$
\lim_{\delta\to 0}\limsup_{N\to\infty} \bb P_{\mu^N} 
\Big[ w_\delta (Y) \ge \varepsilon \Big]\; =\; 0
$$
for every $\varepsilon >0$. Here $w_\delta (Y)$ stands for the 
uniform modulus of continuity defined by
$$
w_\delta (Y)\; =\; \sup_{\substack{ |s-t|\le \delta
\\ \scriptstyle 0\le s, t\le T}} \Vert Y_t -Y_s \Vert_{-k}\; .
$$

We start with a key estimate.  Recall the definition of the
martingales $M_t^{1,N} (H)$, $M_t^{2,N}(H)$ defined by \eqref{eq:01}. 

\begin{lemma}
\label{s21}
Fix a sequence of probability measures $\{\mu^N :\, N\ge 1\}$
satisfying \eqref{f21}, \eqref{f22}.  There exists a finite constant
$C_1$, depending only on $C_0$, such that for every $j\ge 1$,
\begin{equation*}
\limsup_{N\to\infty} \bb E_{\mu^N} \Big[ \sup_{0\le t\le T} 
Y_t(e_j)^2 \Big] \; \le \;  C_1 \, j^4 \, (1+T)^2 \; .
\end{equation*}
\end{lemma}

\begin{proof}
Recall \eqref{eq:01} and write $Y_t^N(e_j)$ as $M_t^{1,N}(e_j) +
Y_0^N(e_j) + \int_0^t Y_s^N (\Delta_N e_j) ds$.  We estimate these three
terms separately.

It follows from \eqref{f21} that $\bb E_{\mu^N} [ Y_0(e_j)^2 ]$ is
bounded by a finite constant $C_1$, uniformly in $N$ and $j$.

Since $M_t^{1,N}(e_j)$ is a martingale, by Doob inequality, $$
\bb
E_{\mu^N} \Big[ \sup_{0\le t\le T} | M_t^{1,N}(e_j) |^2 \Big]\; \le \;
4\, \bb E_{\mu^N} \Big[ \, | M_T^{1,N}(e_j) |^2 \Big]\; .  $$
By
definition of the martingale $M_T^{2,N}(e_j)$ and by \eqref{eq:06},
the right hand side is equal to 
$$
4 \, \bb E_{\mu^N} \Big[ \int_0^T
ds\, \frac 1N \sum_{x=1}^{N-2} [\eta_s(x+1) - \eta_s(x)]^2 \big\{
(\nabla_N e_j)(x/N) \big\}^2 \Big] \;+\; O(N^{-1}) \; .
$$
By Lemma \ref{rl}, as $N\uparrow\infty$, this expression converges
to $8 \int_0^T dt \int_0^1 \chi( \rho(t,u)) (\nabla e_j(u))^2\; du$
which is bounded by $C_1 T j^2$.

Finally, by definition of $e_j$ and by Schwarz inequality,
\begin{equation}
\label{f28}
\bb E_{\mu^N}\Big[ \sup_{0\le t\le T} 
\Big( \int_0^t ds\, Y_s(\Delta e_j) \Big)^2 \Big]  
\;\le\; C_1 \, j^4 \, T \, \bb E_{\mu^N}\Big[ \int_0^T ds
\, Y_s (e_j)^2 \Big] 
\end{equation}
for some finite constant $C_1$. The previous expectation can be
rewritten as
\begin{eqnarray*}
\!\!\!\!\!\!\!\!\!\!\!\!\! &&
\int_0^T ds \, \frac 1N \sum_{x\in\Lambda_N} e_j(x/N)^2 
\rho^N_s(x) [1-\rho^N_s(x)] \\
\!\!\!\!\!\!\!\!\!\!\!\!\! && \quad
+\; \int_0^T ds \, \frac 1N \sum_{x,y\in\Lambda_N} e_j(x/N)
e_j(y/N) \varphi^N_s(x,y) \;, 
\end{eqnarray*}
where
\begin{equation}
\label{f27}
\varphi^N_t(x,y) \;=\; E_{\mu^N}[\eta_t(x); \eta_t(y)]\;.
\end{equation}
By Proposition \ref{s34}, 
\begin{equation}
\label{f29}
\sup_{N\ge 2} \sup_{t\ge 0} N |\varphi^N_t(x,y)| \;\le\; C_1
\end{equation}
for $C_1 = C_0 + (1/2)C_0^2$. Hence, \eqref{f28} is bounded above by
$C_1 T^2 j^4$, which concludes the proof of the lemma.
\end{proof}

\begin{corollary}
\label{s2}
For $k>5/2$,
\begin{eqnarray*}
\!\!\!\!\!\!\!\!\!\!\! &&
\text{(a)}\quad \limsup_{N\to\infty} 
\bb E_{\mu^N}\Big[ \sup_{0\le t\le T}
\Vert Y_t\Vert^2_{-k} \Big]\; <\; \infty\;, \hfill \\
\!\!\!\!\!\!\!\!\!\!\! &&
\text{(b)}\quad \lim_{n\to\infty} \limsup_{N\to\infty} 
\bb E_{\mu^N}\Big[ \sup_{0\le t\le T} \sum_{j\ge n} 
Y_t(e_j)^2 j^{-2k} \Big]\; =\; 0\;.
\end{eqnarray*}
\end{corollary}

The proof of this result is similar to the one of Corollary XI.3.5 in
\cite{kl} and therefore omitted.

In view of Lemma \ref{s21} and part (b) of Corollary \ref{s2}, in order
to prove that the sequence $Q_N$ is tight, we only have to show that
$$
\lim_{\delta\to 0} \limsup_{N\to\infty} \bb P_{\mu^N} 
\Big[ \sup_{\substack{ 0\le|s-t|\le \delta\\ 
0\le s,t\le T}} \{ Y_t(e_j) - Y_s (e_j) \}^2  >\varepsilon\Big]\; =\; 0
$$
for every $j\ge 1$ and $\varepsilon>0$.  Fix $j\ge 1$ and recall
the definition of the martingale $M^{1,N}_t(e_j)$. Since $Y^N_t(e_j)
= Y^N_0(e_j) + M^{1,N}_t(e_j) + \int_0^t \Gamma^N_{s} (e_j) ds$,
the previous statement follows from the next two claims: For every
function $G$ in $C^2_0([0,1])$ and every $\varepsilon>0$,
\begin{eqnarray*}
\!\!\!\!\!\!\!\!\!\!\! &&
\lim_{\delta\to 0}\limsup_{N\to\infty} \bb P_{\mu^N} \Big[
\sup_{\substack{ 0\le s, t\le T \\ \scriptstyle |t-s|\le \delta}}
|M^{1,N}_t(G) - M^{1,N}_s(G)| >\varepsilon \Big] \; =\; 0\; , \\
\!\!\!\!\!\!\!\!\!\!\! &&
\lim_{\delta\to 0}\limsup_{N\to\infty} \bb P_{\mu^N} \Big[
\sup_{\substack{ 0\le s, t\le T \\ \scriptstyle |t-s|\le \delta}}
\Big\vert \int_s^t  dr\, Y^N_r(\Delta_N G) \Big\vert >\varepsilon \Big] \; 
=\; 0\; .
\end{eqnarray*}
The derivation of these estimates is similar to the proofs of Lemmata
XI.3.7 and XI.3.8 in \cite{kl} if one keeps in mind the arguments
presented in the proof of Lemma \ref{s1} and the bound \eqref{f29} on
the two point correlation function $\varphi^N_t$ given by \eqref{f27}.

\section{Semidiscrete heat equation}
\label{sec3}

We prove in this section a bound on the two point correlation function
$\varphi^N_t(x,y)$ introduced in \eqref{f27}. Throughout this section,
$N\ge 2$ is fixed.

For the square of points $C=\{0,\cdots,N\}^2$, consider the subsets
$V=\{(x,y)\in C:\;0<x<y<N\}$ and its boundary $\partial V=\{(x,y)\in
C:\;x=0 \textrm{ or }y=N\}$. Let $\mathcal{M}=\{f:V\cup \partial
V\mapsto \bb R:\; f\big|_{\partial V}=0\}$ and denote by $\Delta_V^N$ the
discrete Laplacian on $\mc M$ defined by
\begin{equation*}
(\Delta_V^N f) (x,y) \;=\; N^2 \big\{
f(x+1,y)+f(x-1,y)+f(x,y-1)+f(x,y+1)-4f(x,y) \big\} 
\end{equation*}
if $|x-y|>1$ and
\begin{equation*}
(\Delta_V^N f) (x,x+1) \;=\; N^2 \big\{
f(x-1,x+1)+f(x,x+2)-2f(x,x+1) \big\} \; .
\end{equation*}
$\Delta_V^N$ corresponds to the generator of a symmetric random walk
on $V\cup \partial V$ which is absorbed on $\partial V$.

We start with an explicit formula for the total time spent by the
random walk on the diagonal, which is expressed by the Green function
or as the solution of the elliptic equation $(- \Delta_V^N \varphi^N)
(x,y) = C \delta_{y=x+1}$, $C$ in $\bb R$. Let $\varphi^N$ be the
solution of
\begin{equation}
\label{eq:05}
\left\{
\begin{array}{l}
\vphantom{\Big\{}
(- \Delta_N \varphi^N)(x,y) = C \, \delta_{y=x+1} 
\quad\text{for $(x,y)$ in $V$}\;, \\
\vphantom{\Big\{}
\varphi^N (x,y) = 0 \quad\text{for $(x,y)$ in $\partial V$}\;.
\end{array}
\right.
\end{equation}
An elementary analysis shows that the unique solution of \eqref{eq:05}
is given by
\begin{equation}
\label{f25}
\varphi^N (x,y)\;=\; \frac {C^2}{N-1} \,
\frac xN \, \Big( 1- \frac yN\Big)\;.
\end{equation}

We turn now to maximum principles for solutions of homogeneous
semidiscrete parabolic equations. Fix a function $\rho^N: \Lambda_N\to
\bb R$ and let $\rho^N_t$ be the solution of
\begin{equation}
\label{f33} 
\left\{
\begin{array}{l}
\partial_s \rho^N_s(x) = \Delta_N \rho_s (x)\; , \;\; x\in\Lambda_N \\
\rho^N_0(x) = \rho^N (x) \; , \;\; x\in\Lambda_N \\ 
\rho^N_s(0) = \alpha \;,\;\; \rho^N_s(N) = \beta \;,\;\; s\ge 0\;.
\end{array}
\right.
\end{equation}

\begin{lemma}
\label{s33}
Let $\rho^N_s$ be the solution of \eqref{f33}. Then,
\begin{equation*}
\sup_{s\ge 0} \max_{0\le x\le N-1} \big | (\nabla_N \rho^N_s) (x) \big
| \;\le\; \max_{0\le x\le N-1} \big | (\nabla_N \rho^N) (x) \big |\; .
\end{equation*}
\end{lemma}

\begin{proof}
Fix $T\ge 0$. Let $\gamma_t(x)=\rho_t(x+1)-\rho_t(x)$ for $0\le x\le
N-1$. Since $\partial_t \gamma_t(x) = (\Delta_N \gamma_t) (x)$ for
$1\le x\le N-2$, by the maximum principle,
\begin{eqnarray*}
\!\!\!\!\!\!\!\!\!\!\!\!\! &&
M \;=\; \max_{0\leq x\leq N-1} \, \sup_{0\leq s\leq T} |\gamma_s(x)| 
\; = \\
\!\!\!\!\!\!\!\!\!\!\!\!\! && \quad
\max \Big \{ \max_{0\leq x\leq N-1} |\rho^N (x+1)-\rho^N (x)|
\, , \, \sup_{0\leq t\leq T} |\rho_t(1)-\alpha| \ ,\, 
\sup_{0\leq t\leq T}|\beta-\rho_t(N-1)| \Big\}\;.
\end{eqnarray*}

We claim that the maximum is attained at $t=0$.  To show this assume,
without loss of generality, that there exists $t_0\in (0,T]$ such that
$M=|\rho_{t_0}(1)-\alpha|$.  By \eqref{f33} with $x=1$ we have,
\begin{equation*}
\partial_s\rho_s(1)=-N^2\rho_s(1)+N^2\alpha+N^2(\rho_s(2)-\rho_s(1))
\end{equation*}
for any $0<s<T$. Thus, multiplying by $e^{sN^2}$, grouping the terms
conveniently and integrating on $[0,t]$ we get that
\begin{equation*}
\rho_t(1)=e^{-tN^2}\rho_0(1)+\alpha(1-e^{-tN^2})+\int_0^t
e^{-(t-s)N^2}N^2(\rho_s(2)-\rho_s(1))ds
\end{equation*}
so that
\begin{equation*}
\rho_t(1)-\alpha=e^{-tN^2}(\rho_0(1)-\alpha)+\int_0^t
e^{-(t-s)N^2}N^2(\rho_s(2)-\rho_s(1))ds\; .
\end{equation*}
Using the assumption made on $|\rho_{t_0}(1)-\alpha|$, we deduce from
this identity that
\begin{equation*}
M \;=\; |\rho_{t_0}(1)-\alpha|
\;\le\; e^{-t_0 N^2}|\rho_0(1)-\alpha| \; + \; (1-e^{-t_0N^2})M\; ,
\end{equation*}
which reduces to $M\leq |\rho_0(1)-\alpha|$. This concludes the proof
of the lemma.
\end{proof}

Fix a function $h: V \to \bb R$ and denote by $f_s$ the solution of
the semidiscrete heat equation
\begin{equation}
\label{f32} 
\left\{
\begin{array}{l}
\partial_s f_s(x,y) = \Delta_V^N f_s(x,y) \;, (x,y)\in V \;, \\ 
f_0(x,y) = h(x,y)\;, \;\; (x,y) \in V \\ 
f_s(x,y) = 0 \;, \;\;  0 \leq s\leq T \;,\;\; (x,y) \in \partial V\; .
\end{array}
\right.
\end{equation}

\begin{lemma}
\label{s32} 
Let $f$ satisfy (\ref{f32}). Then, the maximum value of $f$ on
$[0,T]\times (V\cup \partial V)$ is attained at a point
$(t_0,x_0,y_0)$ such that $t_0=0$ or $(x_0,y_0) \in \partial V$.
\end{lemma}

\begin{proof}
The proof is the same as that of the maximum principle for the usual
heat equation. It uses that if the maximum is attained at an interior
point $(t_0,x_0,y_0)$ then $\partial_{t_0}f_{t_0}(x_0,y_0) \geq 0$ and
$\Delta_V^N f_{t_0}(x_0,y_0)\leq 0$.
\end{proof}

Fix a function $h: V \to \bb R$ and a function $g: \bb R_+\times V
\to\bb R$. Consider the following nonhomogeneous parabolic equation,
\begin{equation}
\label{f31}
\left\{
\begin{array}{l}
\partial_s\varphi_s(x,y) = \Delta_V^N\varphi_s(x,y) + g_s(x,y)\; ,
\;\; (x,y)\in V \;, \\ 
\varphi_0(x,y) = h(x,y)\; , \;\; (x,y)\in V \;, \\
\varphi_s(x,y) = 0\;, \;\;  0\leq s\leq T \;, \;\; (x,y) \in \partial
V\; ,
\end{array}
\right.
\end{equation}

Denote by $\|\cdot \|_{l^{\infty}(V)}$ the sup norm:
$\|h\|_{l^{\infty}(V)} = \max_{z\in V} |h(z)|$.

\begin{lemma}
\label{s31}
Fix $T>0$ and assume that the function $g_s$ is supported on the line
$y=x+1$. Then,
\begin{equation*}
\sup_{0\leq t\leq T} \|\varphi_{t} \|_{l^{\infty}(V)} 
\;\leq\; \| h \|_{l^{\infty}(V)}  
\; +\; \frac{1}{4(N-1)} \sup_{0\leq t\leq T} 
\max_{1\leq x\leq N-2}|g_t(x,x+1)| \; . 
\end{equation*}
\end{lemma}

\begin{proof}
It is not difficult to see that $\Delta_V^N$ is a symmetric, negative
operator on $\mathcal{M}=\{f:V\cup \partial V\mapsto \bb R:\;
f\big|_{\partial V}=0\}$. It generates in particular a strongly
continuous semigroup $\{ e^{s\Delta_V^N} \, : s\geq 0\}$ on
$\mathcal{M}$ and the solution $\varphi_s$ of \eqref{f31} can be
written in the form
\begin{equation*}
\varphi_s \;=\; e^{s\Delta_V^N}h \;+\; 
\int_0^s e^{(s-t)\Delta_V^N} g_t \, dt \; .
\end{equation*}
Since $e^{s\Delta_V^N}h$ is the solution of \eqref{f32}, by the
maximum principle stated in Lemma \ref{s32}, 
\begin{equation*}
\|e^{s\Delta_V^N}h\|_{l^{\infty}(V)} \;\leq\; \|h\|_{l^{\infty}(V)}
\end{equation*}
and $e^{s\Delta_V^N} (x,y) \ge 0$. On the other hand, since $g_t$ is
supported on the line $y=x+1$,
\begin{eqnarray*}
\!\!\!\!\!\!\!\!\!\!\!\!\!\! &&
\Big| \int_0^s \{e^{(s-t)\Delta_V^N} g_t\} (x,y) \, dt \Big| 
\;\le\;  \int_0^s \sum_{x'=1}^{N-2} e^{(s-t)\Delta_V^N}(x,y,x',x'+1)\, 
|g_t(x',x'+1)| \, dt \\ 
\!\!\!\!\!\!\!\!\!\!\!\!\!\! && \qquad
\leq\; \int_0^s \sum_{x'=1}^{N-2} e^{t \Delta_V^N} (x,y,x',x'+1) \, dt 
\sup_{0\leq t\leq T}\max_{1\leq x'\leq N-2}|g_t(x',x'+1)|
\end{eqnarray*}
provided $s\le T$. To conclude the proof of the lemma it remains to
show that the integral is bounded by $(4 [N-1])^{-1}$, uniformly in
$x$, $y$ and $s$.  This integral is bounded above by
\begin{equation}
v_N(x,y) \;=\; \int_0^{\infty} \sum_{x'=1}^{N-2} e^{t \Delta_V^N} 
(x,y,x',x'+1)\, dt 
\end{equation}
which satisfies the equation $\Delta_V^N \, v_N(x,y) =
-\delta_{y=x+1}$. By \eqref{f25},
\begin{equation*}
v_N(x,y) \;=\; \frac{1}{N-1}\frac{x}{N}\left(1-\frac{y}{N}\right)
\end{equation*}
and this expression is less than or equal to $(4 [N-1])^{-1}$ because
$x<y$. This concludes the proof of the lemma.
\end{proof}

We are now in a position to state the main result of this section.
\begin{proposition}
\label{s34}
Consider a probability measure $\mu^N$ on $\{0,1\}^{\Lambda_N}$
satisfying \eqref{f21}. Denote by $\varphi^N_t$ the two point
correlation function:
\begin{equation*}
\varphi^N_t (x,y) \; =\; E_{\mu^N}\Big[  
\{\eta_t (x) - \rho_t^N (x)\} \, \{\eta_t (y) -
\rho_t^N (y)\}  \Big]\;, 
\end{equation*}
where $\rho_t^N$ is the solution of \eqref{f33}. Then,
\begin{equation*}
\sup_{t\ge 0} \|\varphi^N_{t} \|_{l^{\infty}(V)} 
\;\leq\; \frac{2C_0 + C_0^2}{2N}
\end{equation*}
for all $N\ge 2$.
\end{proposition}

\begin{proof}
A simple computation shows that $\varphi^N_t (x,y)$ is the solution of
\eqref{f31} with $h=E_{\mu^N} [\eta (x) ; \eta (y)]$ and $g_t(x,y) = -
(\nabla_N \rho^N_t)(x)^2 \delta_{y=x+1}$. By Lemma \ref{s33} and
assumption \eqref{f21}, $g_t$, which is supported on the line $y=x+1$,
is absolutely bounded by $C_0^2$. Therefore, by Lemma \ref{s31},
$N \sup_{t\ge 0} \|\varphi^N_{t} \|_{l^{\infty}(V)}$ is less than or
equal to $C_0 + C_0^2/2$ since $N\ge 2$. This concludes the proof of
the lemma.
\end{proof}

It remains to check that stationary state $\nu^N_{\alpha, \beta}$
satisfies the assumptions of the previous proposition.  Recall the
definition of $\rho^N(x)$ given just before \eqref{eq:04} and denote
by $\varphi^N (x,y)$, $1\le x < y \le N-1$ the two point correlation
function:
\begin{equation}
\label{eq:03}
\varphi^N (x,y) \; =\; E_{\nu^N_{\alpha, \beta}}\Big[  
\{\eta (x) - \rho^N (x)\} \, \{\eta (y) -
\rho^N (y)\}  \Big]\;. 
\end{equation}
Since $ E_{\nu^N_{\alpha, \beta}} [ L_N \{\eta (x) -
\rho^N (x)\} \, \{\eta (y) - \rho^N (y)\} ]=0$ for all $1\le x < y \le
N-1$, we obtain that $\varphi^N (x,y)$ is the solution of the discrete
differential equation \eqref{eq:05} with $C= \beta -
\alpha$. Therefore,
\begin{equation*}
\varphi^N (x,y)\;=\; \frac {(\beta-\alpha)^2}{N-1} \,
\frac xN \, \Big( 1- \frac yN\Big)\;.
\end{equation*}
In particular, \eqref{f21} is satisfied and we may apply Proposition
\ref{s34}.

\end{document}